\documentclass[12pt]{amsart}
\usepackage{amsfonts}
\usepackage{amsmath,amssymb}

\newtheorem{theorem}{Theorem}[section]
\newtheorem{lemma}[theorem]{Lemma}
\newtheorem{claim}[theorem]{Claim}

\newtheorem{corollary}[theorem]{Corollary}

\theoremstyle{definition}

\newtheorem{remark}[theorem]{Remark}

\DeclareMathOperator{\Vol}{Vol} 
 
\DeclareMathOperator{\Hess}{Hess} 
 
\DeclareMathOperator{\diam}{diam} 
\DeclareMathOperator{\Diff}{Diff}

\begin{document}

\title[Degeneration of shrinking Ricci solitons]
{Degeneration of shrinking Ricci solitons}
\author{Zhenlei Zhang}
\address{Department of Mathematics, Capital Normal University, Beijing, 100037, P. R. China}
\email{zhleigo@yahoo.com.cn}
\thanks{The author was supported by the Capital Normal University}

\maketitle


\begin{abstract}
Let $(Y,d)$ be a Gromov-Hausdorff limit of closed shrinking Ricci
solitons with uniformly upper bounded diameter and lower bounded
volume. We prove that off a closed subset of codimension at least 2,
$Y$ is a smooth manifold satisfying a shrinking Ricci soliton
equation.
\end{abstract}

\section{Introduction}

Let $(Y,d)$ be a metric space obtained as a Gromov-Hausdorff limit
of a sequence of complete $n$-dimensional Riemannian manifolds
$(M_k,g_k)$. When the Ricci curvatures of $M_k$ are uniformly
bounded below, the work of Cheeger and Colding \cite{ChCo1, ChCo2,
ChCo3} gives a beautiful description for the structure of $Y$. In
particular they proved in the noncollapsing case that $Y$ is
bi-H\"{o}lder equivalent to connected smooth Riemannian manifold
outside of a singular set of codimension $\geq2$. When the Ricci
curvatures of $M_k$ are absolutely bounded, they also proved that
the singular set of $Y$ is closed.

As a generalization of notion, Bakry-\'{E}mery Ricci curvature plays
an important role in the theory of smooth measure metric space, cf.
\cite{WeWy2} and the references therein. In view of Cheeger and
Colding's work, one naturally asks what will happen to $Y$ when it
is obtained as a limit of Riemannian manifolds with bounded
Bakry-\'{E}mery Ricci curvature. In this paper, we are going to
prove a partial result to this question, for the special case when
$M_k$ are closed shrinking Ricci solitons. The argument here also
applies to the noncompact Ricci solitons, either for shrinking,
steady or expanding case.

The Ricci solitons are manifolds with constant Bakry-\'{E}mery Ricci
curvature. Namely, a \textit{Ricci soliton} is a Riemannian manifold
$(M,g)$ satisfying
\begin{equation}
Ric+\Hess(f)=\epsilon g
\end{equation}
for some $\epsilon\in\mathbb{R}$ and $f\in C^\infty(M;\mathbb{R})$.
The associated function $f$ is called the \textit{potential
function}. We assume $\epsilon=+\frac{1}{2},-\frac{1}{2}$ or $0$
after a normalization, in each case the Ricci soliton is called a
\textit{shrinking, steady} or \textit{expanding} one respectively.
Ricci solitons are also known as the self-similar solutions to the
Ricci flow, where the metric varies via rescalings and the
diffeomorphic transformations, cf. \cite{CK}.

Ricci solitons are generalizations of Einstein manifolds. Suppose
$M_k$ as in the later case, Cheeger and Colding proved in
\cite{ChCo1} that, under the noncollapsing hypothesis, $Y$ is
Einstein off the singular set and the convergence is smooth off the
singular set. We will prove the following Ricci soliton analogy of
Cheeger and Colding's result in the present paper.


\begin{theorem}\label{t1}
Let $(M_k,g_k)$ be a sequence of n-dimensional closed shrinking
Ricci solitons satisfying {\rm(1)} with $\epsilon=\frac{1}{2}$.
Suppose that
\begin{equation}\label{volume noncollapsing}
\Vol_{g_k}(M_k)\geq v,
\end{equation}
\begin{equation}\label{diameter bound}
\diam(M_k,g_k)\leq D
\end{equation}
uniformly hold for some $0<v,D<\infty$. Then passing a subsequence
if necessary, the manifolds $(M_k,g_k)$ converge in the
Gromov-Hausdorff sense to a length metric space $(Y,d_\infty)$.

The singular set $\mathcal{S}$ is closed in $Y$ and has Hausdorff
codimension at least 2; the regular set
$\mathcal{R}=Y\backslash\mathcal{S}$ is connected. On $\mathcal{R}$,
the metric is induced from a smooth Riemannian metric which
satisfies a shrinking Ricci soliton equation and the convergence
takes place smoothly.
\end{theorem}

By definition, cf. \cite{ChCo1}, a point $y$ belongs to
$\mathcal{S}$ iff there is a tangent cone at $y$ which is not
isometric to the $n$-Euclidean space.

The proof of the theorem proceeds as follows: We first show the
theorem up to a conformal change of the shrinking Ricci solitons, by
using Cheeger and Colding's theorem on degeneration for Ricci
curvature \cite{ChCo1}; then we show that the degeneration property
for the shrinking Ricci solitons does not depend on the conformal
changes. The point is that the approximation maps associated to the
Gromov-Hausdorff convergence do not depend on the conformal change
of the metrics.

When the curvatures have bounded $L^{n/2}$-norms, the singularities
are finite and of orbifold type; see \cite{We}, \cite{ChWa} or
\cite{Zhang2}. Applying Theorem \ref{t1} together with Theorem 2.6
of \cite{An2}, we give a more direct proof of this result; see \S
3.6. We remark that the first orbifold type convergence theorem for
Ricci solitons was given by Cao and Sesum in the K\"{a}hler case
\cite{CaSe}, under additional assumption of lower bounded Ricci
curvature. Later, X. Zhang proved the analogy for real case
\cite{Zh}. The orbifold type compactness theorem is a generalization
of Einstein manifolds, see \cite{An1}, \cite{BaKaNa} and \cite{Ti1}.


In the K\"{a}hler category, one may hope to know more about the
structure of $Y$, as in the Ricci curvature case (cf.
\cite{ChCoTi}). Unfortunately, our method in this paper can not
applied to the K\"{a}hler manifolds.

However, by a different argument, we can improve the orbifold
compactness theorem mentioned above as in the following theorem for
K\"{a}hler Ricci solitons. The result is a shrinking Ricci soliton
version of Tian's compactness theorem for K\"{a}hler Einstein
manifolds \cite{Ti2}. It can also be seen as a sharpening of Cao and
Sesum's orbifold compactness theorem for shrinking K\"{a}hler Ricci
solitons \cite{CaSe} of dimension at least 3.
\begin{theorem}\label{t2}
Let $(M_k,g_k)$ be a sequence of K\"{a}hler Ricci solitons with
positive first Chern class, of dimension
$n=\dim_{\mathbb{C}}M_k\geq3$. Suppose {\rm(\ref{volume
noncollapsing})} and {\rm(\ref{diameter bound})} hold and that
\begin{equation}\label{integral curvature}
\int_{M_k}|Rm(g_k)|^ndv_{g_k}\leq C
\end{equation}
for some $C$ independent of $k$, then passing a subsequence,
$(M_k,g_k)$ converge smoothly to another K\"{a}hler Ricci soliton of
positive first Chern class.
\end{theorem}

\begin{remark}
Under the hypothesis of bounded $L^{n/2}$-norm of curvature tensor,
the condition (\ref{volume noncollapsing}) together with
(\ref{diameter bound}) is equivalent to the boundedness of
Perelman's $\nu$ functional; see \cite{CaSe} or \cite{Zhang1} for a
lower bound of $\nu$ in terms of (\ref{volume noncollapsing}) and
(\ref{diameter bound}), and see \cite{Zhang2} for a lower bound of
volume and upper bound of diameter via $\nu$ and (\ref{integral
curvature}).
\end{remark}

At the end of the introduction, we give a remark on the geometry of
Bakry-\'{E}mery Ricci curvature.

\begin{remark}
Bakry-\'{E}mery Ricci curvature shares many similar properties as
the Ricci curvature, especially when the potential function admits a
$C^0$ bound. The key point is that, when Laplacian is replaced by a
modified elliptic operator, the Laplace comparison theorem for
distance function remains true up to a perturbation which is
controlled by the $C^0$ norm of the potential, cf. \cite{FLZ}. In
addition, the Bochner formula for Bakry-\'{E}mery Ricci curvature
remains valid \cite{BaEm}; see also \cite{FLZ, WeWy2}. With these
two things in hand, many results for Ricci curvature can be extended
to Bakry-\'{E}mery Ricci curvature, such as the splitting theorem
\cite{FLZ, WeWy1}.
\end{remark}

The paper is organized as follows: In \S 2, we recall some
preliminaries that will be used; in \S 3, we give a proof of Theorem
\ref{t1}. In \S 4, we consider the K\"{a}hler case and present a
proof of Theorem \ref{t2}.

\vskip 2mm

\noindent {\bf Acknowledgement:} The author would like to thank
professor X.D. Li for asking him to pay attention to the question on
the extension of Cheeger and Colding's theory to Bakry-\'{E}mery
Ricci curvature. Thanks also goes to professor Fuquan Fang for his
constant help and encouragement. The author also thanks Yuguang
Zhang for his comments on the proof of Theorem \ref{t2}.

\section{Preliminaries}

\subsection{Pseudolocality theorem for Ricci flow}

We state another version of Perelman's pseudolocality theorem which
is proved in \cite{FZZ}. The difference from Perelman's theorem
\cite{Pe} is that here we use the local almost Euclidean volume
growth instead of almost Euclidean isoperimetric estimate.
\begin{theorem}[\cite{FZZ}]\label{Pseudolocality}
There exist universal constants $\delta_{0},\epsilon_{0}>0$ with the
following property. Let $g(t),t\in[0,(\epsilon_{0}r_{0})^{2}],$ be a
solution to the Ricci flow
\begin{equation}
\frac{\partial}{\partial t}g(t)=-2Ric(g(t))
\end{equation}
on a closed $n$-manifold $M$ and $x_{0}\in M$ be a point. If the
scalar curvature
$$R(x,t)\geq-r_{0}^{-2}\mbox{ whenever }{\rm{dist}}_{g(t)}(x_{0},x)\leq
r_{0},$$ and the volume
$${\rm{Vol}}_{g(t)}(B_{g(t)}(x,r))\geq(1-\delta_{0}){\rm{Vol}}(B(r))\mbox{ for all }B_{g(t)}(x,r)\subset B_{g(t)}(x_{0},r_{0}),$$ where
$B(r)$ denotes a ball of radius $r$ in the $n$-Euclidean space and
${\rm{Vol}}(B(r))$ denotes its Euclidean volume, then the Riemannian
curvature tensor satisfies $$|Rm|_{g(t)}(x,t)\leq t^{-1},\mbox{
whenever }{\rm{dist}}_{g(t)}(x_{0},x)<\epsilon_{0}r_{0}.$$

In particular, $|Rm|_{g(t)}(x_{0},t)\leq t^{-1}$ for all time
$t\in(0,(\epsilon_{0}r_{0})^{2}]$.
\end{theorem}

We use the pseudolocality theorem to prove the smooth convergence on
the regular part of the limit space.

\subsection{One technical lemma}

Let $(M,g)$ be a shrinking Ricci soliton satisfying (\ref{volume
noncollapsing}) and (\ref{diameter bound}). Then the potential
function has a uniform $C^1$ norm bound; see (\ref{e35}) in \S 3.1.
Using the mean curvature and relative volume comparison theorems for
Bakry-\'{E}mery Ricci curvature, cf. \cite{WeWy1}, one easily
derives the following technical lemma.
\begin{lemma}\label{l21}
For all $\epsilon>0$, there exists $c=c(n,v,D,\epsilon)>0$ such that
the following holds. Let $E\subset M$ be a submanifold with smooth
boundary and $x_1,x_2\in M$ be two points such that
\begin{equation}
B(x_1,\epsilon)\cup B(x_2,\epsilon)\subset M\backslash E.
\end{equation}
If every minimal geodesic from $x_1$ to point in $B(x_2,\epsilon)$
intersects $E$, then
\begin{equation}
\Vol_g(\partial E)\geq c.
\end{equation}
\end{lemma}
\begin{proof}
The proof follows directly from the argument in Page 525 of
\cite{Gr}; see also \cite{ChCo2}. We omit the details here.
\end{proof}

This lemma is applied to prove that on the regular part of the limit
space, the intrinsic metric coincides with the extrinsic metric; see
also \cite{ChCo1}.

\subsection{Isoperimetric inequality}

Let $(M,g)$ be an $n$-dimensional Riemannian manifold. The global
isoperimetric constant is defined to be the minimal constant
$C_I(M,g)$ such that
\begin{equation}\label{iso}
\min(\Vol(\Omega)^{n-1},\Vol(M\backslash\Omega)^{n-1})\leq
C_I(M,g)\Vol(\partial\Omega)^{n}
\end{equation}
for all domain $\Omega\subset M$ with smooth boundary. It is well
known that $C_I(M,g)$ is equivalent to the Sobolev constant, which
is defined to be the minimal constant $C_S(M,g)$ such that
\begin{equation}
\inf_{a\in\mathbb{R}}(\int_M(\psi-a)^{\frac{n}{n-1}}dv)^{\frac{n-1}{n}}\leq
C_S(M,g)\int_M|\nabla \phi|
\end{equation}
for all smooth function $\phi$. It is a direct calculation that for
equivalent metrics $g$ and $\tilde{g}$ with relation
$$C^{-1}g\leq\tilde{g}\leq Cg,$$
the Sobolev constants satisfies
\begin{equation}
C^{1-n}\cdot C_S(M,g)\leq C_S(M,\tilde{g})\leq C^{n-1}\cdot
C_S(M,g).
\end{equation}
The isoperimetric constants for $g$ and $\tilde{g}$ are also
equivalent to each other.

If $(M,g)$ is a closed shrinking Ricci soliton which satisfies
(\ref{volume noncollapsing}) and (\ref{diameter bound}), then the
metric $\tilde{g}=e^{-\frac{2}{n-2}f}g$ satisfies
$$C_2^{-1}g\leq\tilde{g}\leq C_2g,\hspace{0.3cm}|Ric(\tilde{g})|_{\tilde{g}}\leq C_3,$$
for some constants $C_2,C_3$ depending on only $n$ and $D$, see
Section 3. Thus the volume $\Vol(\tilde{g})\geq C_2^{-n/2}v$ and the
diameter $\diam(\tilde{g})\leq C_2^{1/2} D$. So, by one of Croke's
classical theorem \cite{Cr}, $C_I(M,\tilde{g})$ as well as
$C_S(M,\tilde{g})$ has a uniform upper bound. Passing to $g$, we
obtain an upper bound
\begin{equation}\label{e20}
C_I(M,g)\leq C_0=C_0(n,v,D).
\end{equation}

\section{Degeneration of Ricci solitons}

Let $(M,g)$ be a closed shrinking Ricci soliton with potential
function $f$:
\begin{equation}\label{soliton}
Ric+\Hess(f)=\frac{1}{2}g.
\end{equation}
Tracing the soliton equation we get
\begin{equation}\label{e31}
R+\triangle f=\frac{n}{2},
\end{equation}
then doing integration yields
\begin{equation}\label{e32}
\frac{1}{\Vol_g(M)}\int_M Rdv=\frac{n}{2}.
\end{equation}
The following equation will also be frequently used, cf. \cite[Eq.
4.13]{ChLuNi}:
\begin{equation}\label{e33}
R+|\nabla f|^2=f+const.,
\end{equation}
Assume that $const.=0$ after a translation.

\subsection{Bound the potential}

We suppose in this section that the diameter satisfies
\begin{equation}
\diam(M,g)\leq D.
\end{equation}
By Ivey \cite{Iv}, a closed shrinking Ricci soliton has positive
scalar curvature. (\ref{e33}) tells us that the infimums of $f$ and
$R$ are attained at the same point. As a consequence, $0<\inf f=\inf
R\leq\frac{n}{2}$. Then using (\ref{e33}) once again one derives the
upper bound of $f$:
\begin{equation}\label{e34}
\sup f\leq\frac{1}{4}(D+\frac{n}{2})^2.
\end{equation}
It follows immediately from (\ref{e33}) that for
$C_1=\frac{3}{4}(D+\frac{n}{2})^2$,
\begin{equation}\label{e35}
\sup f+\sup|\nabla f|^2+\sup R\leq C_1.
\end{equation}

\subsection{Conformal change of the soliton metric}

Define a conformal metric $\tilde{g}=e^{-\frac{2}{n-2}f}g$. It
follows an easy computation that
\begin{equation}\label{e36}
C_2^{-1}g(x)\leq\tilde{g}(x)\leq g(x),
\end{equation}
where $C_2=e^{\frac{(D+\frac{n}{2})^2}{2(n-2)}}$ depending on $n,D$.
Let $\widetilde{Ric}$ denote the Ricci tensor of $\tilde{g}$, then,
together with using equations (\ref{e31}) and (\ref{e33}), we get
that, cf. \cite{Be},
\begin{eqnarray}
\widetilde{Ric}&=&Ric+\Hess(f)+\frac{1}{n-2}d
f\otimes d f+\frac{1}{n-2}(\triangle f-|\nabla f|^{2})g\nonumber\\
&=&\frac{1}{2}g+\frac{1}{n-2}d f\otimes d f
+\frac{1}{n-2}(\frac{n}{2}-R(g)-|\nabla f|^{2})g\nonumber\\
&=&\frac{1}{n-2}(d f\otimes d f+(n-1-f)g).\nonumber
\end{eqnarray}

So we can bound the Ricci curvature after a conformal change:

\begin{lemma}\label{l31}
There exists $C_3$ depending only on $n$ and $D$ such that the Ricci
curvature of $\tilde{g}$ satisfies the bound
\begin{equation}\label{e37}
|\widetilde{Ric}|_{\tilde{g}}\leq C_3.
\end{equation}
\end{lemma}
\begin{proof}
By a straightforward computation,
\begin{eqnarray}
|\widetilde{Ric}|_{\tilde{g}}&\leq&\frac{1}{n-2}(|df\otimes
df|_{\tilde{g}}+|n-1-f|\cdot|g|_{\tilde{g}})\nonumber\\
&\leq&\frac{C_2}{n-2}(|\nabla f|_g^2+n|n-1-f|).\nonumber
\end{eqnarray}
By (\ref{e35}), one can choose
$C_3=\frac{C_2}{n-2}(n(n-1)+(n+1)C_1)$.
\end{proof}

\subsection{Bound the volume ratio}

Suppose further that
\begin{equation}
\Vol_g(M)\geq v.
\end{equation}
Then $\Vol_{\tilde{g}}(M)\geq C_2^{-n/2}\Vol_{g}(M)\geq
C_2^{-n/2}v$. By relative volume comparison theorem, together with
using that $\tilde{g}$ and $g$ are uniformly equivalent, there
exists a positive constant $\kappa=\kappa(n,v,D)$ such that
\begin{equation}\label{e38}
\kappa r^n\leq\Vol_g(B(r))\leq\kappa^{-1}r^n
\end{equation}
for all metric ball $B(r)$ of radius $r\leq D$ in $M$.

\subsection{Convergence modulo conformal changes}

Let $(M_k,g_k)$ be a sequence of n-dimensional closed shrinking
Ricci solitons with potential functions $f_k$. Suppose that
$(M_k,g_k)$ satisfies {\rm (\ref{volume noncollapsing})} and {\rm
(\ref{diameter bound})}. Let $\tilde{g}_k=e^{-\frac{2}{n-2}f_k}g_k$.
Applying Cheeger-Colding's theorem, cf. \cite[Thm. 7.2]{ChCo1},
gives the following
\begin{theorem}
Passing a subsequence if necessary, the Riemannian manifolds
$(M_k,\tilde{g}_k)$ converge in the Gromov-Hausdorff sense to a
compact n-dimensional length metric space $(Y,\tilde{d}_\infty)$.

The singular set $\mathcal{S}$ is a closed subset which has
Hausdorff codimension at least 2. On $Y\backslash\mathcal{S}$, the
metric is induced by a $C^{1,\alpha}$ Riemannian metric for all
$\alpha<1$. Furthermore, the convergence is in the $C^{1,\alpha}$
sense on $Y\backslash\mathcal{S}$.
\end{theorem}

Based on this theorem, we can finally give a

\subsection{Proof of Theorem \ref{t1}}

Let $(M_k,g_k)$ be a sequence of n-dimensional closed shrinking
Ricci solitons satisfying {\rm(\ref{soliton})} with potential
functions $f_k$. Suppose further that $(M_k,g_k)$ satisfies {\rm
(\ref{volume noncollapsing})} and {\rm (\ref{diameter bound})}.

Let $\tilde{g}_k=e^{-\frac{2}{n-2}f_k}g_k$ and
$(Y,\tilde{d}_\infty)$ be the Gromov-Hausdoff limit of
$(M_k,\tilde{g}_k)$. Let $\mathcal{R}$ be the regular part of $Y$
and $\mathcal{S}=M\backslash\mathcal{R}$ be the singular set.
$\mathcal{R}$ is a $C^{1,\alpha}$ manifold on which the convergence
is in the $C^{1,\alpha}$ sense. Let $\tilde{g}_\infty$ be the
Riemannian metric defined on $\mathcal{R}$.


Let $\{K_k\}_{k=1}^\infty$ be any given exhaustion of $\mathcal{R}$
by compact subsets, such that $K_k\subset K_{k+1}$ for all $k$ and
$\cup_{k=1}^\infty K_k=\mathcal{R}$. Then, by definition of
$C^{1,\alpha}$ convergence, there exists a sequence of smooth
embeddings
$$\psi_k:K_k\rightarrow M_k$$
such that
$\psi_k^*\tilde{g}_k\stackrel{C^{1,\alpha}}{\longrightarrow}
\tilde{g}_\infty$ as $k\rightarrow\infty$ on any compact sets of
$\mathcal{R}$.

\begin{claim}\label{c31}
Passing a subsequence if necessary, $\psi_k^*g_k$ converges in the
$C^{\alpha}$ sense to a $C^{\alpha}$ metric, say $g_\infty$, on
$\mathcal{R}$.
\end{claim}
\begin{proof}
It follows from that $\|f_k\|_{C^1}$ is uniformly bounded on $M_k$,
and that the metrics
$$\psi_k^*g_k=\psi_k^*(e^{\frac{2}{n-2}f_k}\tilde{g}_k)=e^{\frac{2}{n-2}f_k\circ\psi_k}\psi_k^*\tilde{g}_k$$
have a uniform $C^1$ bound on any compact subset of $\mathcal{R}$.
\end{proof}

In view of (\ref{e36}), the metrics $g_\infty$ is uniformly
equivalent to $\tilde{g}_\infty$ on $\mathcal{R}$:
\begin{equation}\label{e310}
C_2^{-1} g_\infty\leq\tilde{g}_\infty\leq g_\infty.
\end{equation}
This concludes that, as intrinsic metric spaces,
$(\mathcal{R},g_\infty)$ and $(\mathcal{R},\tilde{g}_\infty)$ have
the same metric completions in the set point of view. In view of
Theorem 7.2 in \cite{ChCo1} and Theorem 3.9 in \cite{ChCo2}, the
intrinsic metric defined on $\mathcal{R}$ in terms of
$\tilde{g}_\infty$ coincides with the extrinsic metric
$\tilde{d}_\infty$, i.e., for all $y_1,y_2\in\mathcal{R}$,
$$\tilde{d}_\infty(y_1,y_2)=\inf\{L_{\tilde{g}_\infty}(\gamma)|\gamma\subset\mathcal{R} \mbox{ is a curve
joining }y_1\mbox{ and }y_2\}.$$ Thus the completion space of
$(\mathcal{R},g_\infty)$, as a set, equals $Y$.

Denote by $d_\infty$ the distance function on $Y$ as the metric
completion of $(\mathcal{R},g_\infty)$, then $\mathcal{R}$ is still
an open set in $(Y,d_\infty)$ by the equivalence of $g_\infty$ and
$\tilde{g}_\infty$. In addition, the distance function $d_\infty$
and $\tilde{d}_\infty$ satisfies
\begin{equation}\label{e311}
C_2^{-1/2}d_\infty(y_1,y_2)\leq \tilde{d}_\infty(y_1,y_2)\leq
d_\infty(y_1,y_2)
\end{equation}
for all $y_1,y_2\in Y$. So the set
$\mathcal{S}=Y\backslash\mathcal{R}$ is closed and has Hausdorff
dimension $\dim(\mathcal{S})\leq n-2$ with respect to $d_\infty$ as
well. We next claim that

\begin{claim}\label{c32}
$(Y,d_\infty)$ is a path metric space.
\end{claim}
\begin{proof}
By definition of metric completion, for any $y_1,y_2\in Y$, there
exist sequences of points $\{y_{1,i},y_{2,i}\}_{i=1}^\infty$ and
minimal geodesics $\gamma_i\subset\mathcal{R}$ connecting $y_{1,i}$
and $y_{2,i}$ such that
$$y_{1,i}\rightarrow y_1,\hspace{0.3cm}y_{2,i}\rightarrow
y_2,\hspace{0.3cm}L_{g_\infty}(\gamma_i)\rightarrow
d_{\infty}(y_1,y_2)$$ as $i\rightarrow\infty$. Obviously the metric
space $(Y,d_\infty)$ is compact, so $\gamma_i$ can be chosen such
that they converge to a minimal geodesic connecting $y_1$ and $y_2$
by Arzela-Ascoli theorem. This implies that $(Y,d_\infty)$ is a path
metric space.
\end{proof}

From now on, we assume that
$(M_k,\tilde{g}_k)\stackrel{d_{GH}}{\longrightarrow}(Y,\tilde{d}_\infty)$;
assume further, as in Claim 3.3, that
$\psi_k^*g_k\stackrel{C^{\alpha}}{\longrightarrow}g_\infty$ on
$\mathcal{R}$. The later convergence is uniform on any compact
subset. Denote by $d_k$ and $\tilde{d}_k$ the distance function
induced by $g_k$ and $\tilde{g}_k$ respectively on $M_k$.

Let $\epsilon>0$ be a fixed constant and $A\subset\mathcal{R}$ be an
$\epsilon$-dense set of $(Y,\tilde{d}_\infty)$, i.e., the
$\epsilon$-neighborhood of $A$ covers $Y$. Obviously $A$ is
$C_1^{1/2}\epsilon$-dense in $(Y,d_\infty)$. The goal for a while is
to show that $\psi_k$ defines a $3\epsilon$-approximation from
$(A,\tilde{d}_\infty)$ to $(M_k,\tilde{g}_k)$ whenever $k$ is large
enough.
\begin{claim}\label{c33}
$\psi_k(A)$ is $3\epsilon$-dense in $(M_k,\tilde{g}_k)$ for $k$
large enough.
\end{claim}
\begin{proof}
Let $K\subset\mathcal{R}$ be a suitable chosen compact submanifold
with smooth boundary such that $A\subset K$ and that
$\Vol_{\tilde{g}_\infty}(\partial K)$ is as small as possible. This
is can be done since the singular set is of codimension at least 2
in $(Y,\tilde{g}_\infty)$. We first show that $\psi_k(K)$ is
$\epsilon$-dense in $(M_k,\tilde{g}_k)$. Let
$\{U_{k,j}\}_{j=1}^{N_k}$ be the components of
$M_k\backslash\psi_k(K)$. Then noting that for $k$ large enough, $K$
is contained in the interior of $K_k$. Thus
$\partial(M_k\backslash\psi_k(K))=\partial\psi_k(K)=\psi_k(\partial
K)$ is a smooth submanifold of $M_k$, so $$\partial
U_{k,j_1}\cap\partial U_{k,j_2}=\emptyset \mbox{ whenever } j_1\neq
j_2.$$ Using the isoperimetric inequality (\ref{iso}) to each domain
$U_{k,j}$, and noticing that $\Vol_{\tilde{g}_\infty}(\partial K)$
is as any small as possible, we obtain that
$$\Vol_{\tilde{g}_k}(U_{k,j})\leq C_I(M_k,\tilde{g}_k)\cdot\Vol_{\tilde{g}_k}(\partial U_{k,j})^{\frac{n}{n-1}}<\kappa\epsilon^n$$
whenever $k$ is large enough. This, together with (\ref{e38}),
implies that $M_k\backslash\psi_k(K)$ contains no balls of radius
$\epsilon$. Thus $\psi_k(K)$ is $\epsilon$-dense in
$(M_k,\tilde{g}_k)$.

On the other hand, since $A$ is $\epsilon$-dense in $K$, we get that
$\psi_k(A)$ is $2\epsilon$-dense in $(\psi_k(K),\tilde{d}_k)$.
Actually, for any $y\in K$, there exists a curve
$\gamma\subset\mathcal{R}$ joining $y$ to a point in $A$ such that
$$L_{\tilde{g}_\infty}(\gamma)\leq\tilde{d}_\infty(y,A)+\epsilon/2<3\epsilon/2.$$
Thus the distance
$$\tilde{d}_k(\psi_k(y),\psi_k(A))\leq L_{\tilde{g}_k}(\psi_k(\gamma))\leq\epsilon/2+L_{\tilde{g}_\infty}(\gamma)<2\epsilon$$
for $k$ large enough, since $\psi_k^*\tilde{g}_k$ converges to
$\tilde{g}_\infty$ in $C^{1,\alpha}$ sense on compact subsets of
$\mathcal{R}$. Then using the fact that $K$ is compact one deduces
the desired result.

Summing up the results obtained completes the proof of the claim.
\end{proof}

Next we show that $\psi_k$ is almost an isometry on $A$.
\begin{claim}\label{c34}
For all $a_1,a_2\in A$ and $k$ large enough,
\begin{equation}\label{e312}
\tilde{d}_k(\psi_k(a_1),\psi_k(a_2))\leq\tilde{d}_{\infty}(a_1,a_2)+2\epsilon.
\end{equation}
\end{claim}
\begin{proof}
Indeed, for $a_1,a_2\in A\subset\mathcal{R}$, there exists a curve
$\gamma\subset\mathcal{R}$ joining $a_1,a_2$ such that
$$L_{\tilde{g}_\infty}(\gamma)\leq\tilde{d}_\infty(a_1,a_2)+\epsilon.$$
On the other hand, the length of curves $\psi_k(\gamma)$ with
respect to $\tilde{g}_k$ is bigger than
$\tilde{d}_k(\psi_k(a_1),\psi_k(a_2))$. Thus, using that
$\psi_k^*\tilde{g}_k$ converges to $\tilde{g}_\infty$ uniformly, we
have for $k$ large enough
$$\tilde{d}_k(\psi_k(a_1),\psi_k(a_2))\leq L_{\tilde{g}_k}(\psi_k(\gamma))\leq L_{\tilde{g}_\infty}(\gamma)+\epsilon
\leq\tilde{d}_\infty(a_1,a_2)+2\epsilon.$$ This proves (\ref{e312}).
\end{proof}

Next we prove the other part:
\begin{claim}\label{c35}
For all $a_1,a_2\in A$ and $k$ large enough.
\begin{equation}\label{e313}
\tilde{d}_k(\psi_k(a_1),\psi_k(a_2))\geq
\tilde{d}_{\infty}(a_1,a_2)-2\epsilon.
\end{equation}
\end{claim}
\begin{proof}
Let $K\subset\mathcal{R}$ be the compact submanifold with smooth
boundary such that $\Vol_{\tilde{g}_\infty}(\partial K)$ is small
enough, as in the proof of Claim \ref{c33}. After a modification of
the $\epsilon$-dense subset $A$, we may assume that $A\subset K$ and
$$\tilde{d}_\infty(A,\partial K)=\min\{\tilde{d}_\infty(a,y)|a\in A,y\in\partial K\}\geq\epsilon/4.$$
Passing to the sequence, we have that
$$\tilde{d}_k(\psi_k(A),\partial\psi_k(K))=\min\{\tilde{d}_k(\psi_k(a),\psi_k(y))|a\in A,\psi_k(y)\in\partial\psi_k(K)\}\geq\epsilon/5,$$
for $k$ large enough. Then for given $a_1,a_2\in A$, if any minimal
geodesic from $\psi_k(a_2)$ to points in
$B_{\tilde{g}_k}(\psi_k(a_1),\epsilon/5)$ intersects
$M_k\backslash\psi_k(K)$, then Lemma \ref{l21} yields that
$$\kappa\cdot(\epsilon/5)^n\leq\Vol_{\tilde{g}_k}(B_{\tilde{g}_k}(\psi_k(a_1),\epsilon/5))\leq
C(n,v,D,\epsilon)\cdot\Vol_{\tilde{g}_k}(\psi_k(\partial K)),$$
which can not happen if $\Vol_{\tilde{g}_\infty}(\partial K)$ is
less than a quantity depending on $n,v,D$ and $\epsilon$. Thus for
any $a_1,a_2\in A$ and $k$ large enough, there exists one minimal
geodesic $\gamma_k$ joining $\psi_k(a_2)$ to a point in
$B_{\tilde{g}_k}(\psi_k(a_1),\epsilon/5)$ which is contained in
$\psi_k(K)$. As $k\rightarrow\infty$, $\gamma_k$ converge to one
minimal geodesic $\gamma_\infty$ contained on $K$ which connects
$a_2$ and one point in $B_{\tilde{g}_\infty}(a_1,\epsilon/5)$. Thus
\begin{eqnarray}
\tilde{d}_k(\psi_k(a_1),\psi_k(a_2))&\geq&
L_{\tilde{g}_k}(\gamma_k)-\epsilon/5 \nonumber\\
&\geq&
L_{\tilde{g}_\infty}(\gamma_\infty)-\epsilon\nonumber\\
&\geq&\tilde{d}_\infty(a_2,\gamma_\infty(1))-\epsilon\nonumber\\
&\geq&\tilde{d}_\infty(a_2,a_1)-2\epsilon\nonumber
\end{eqnarray}
whenever $k$ is large enough. This proves the claim.

We mention that the key point here is that, the submanifold $K$ can
be chosen such that its boundary has volume as small as possible.
\end{proof}

Passing to the the metrics $g_k$, $\psi_k(A)$ is
$3C_1^{1/2}\epsilon$-dense in $M_k$. Applying the same argument as
in Claim \ref{c34} and \ref{c35}, by the $C^\alpha$ convergence on
$\mathcal{R}$, one can show that
$$|d_k(\psi_k(a_1),\psi_k(a_2))-d_\infty(a_1,a_2)|\leq2\epsilon$$
for all $a_1,a_2\in A$ and $k$ large enough. This means that
$\psi_k$ defines a $3C_1^{1/2}\epsilon$-approximation from
$(A,d_\infty)$ to $(M_k,d_k)$. As a consequence, noting that $A$ is
$C_1^{1/2}\epsilon$-dense in $(Y,d_\infty)$ and then letting
$\epsilon\rightarrow0$, we finally obtain
\begin{claim}\label{c36}
$(M_k,g_k)$ converges to $(Y,d_\infty)$ in the Gromov-Hausdorff
sense.
\end{claim}

From the arguments above, the maps $\psi_k$ define the approximation
for both Gromov-Hausdorff convergence
$(M_k,g_k)\stackrel{d_{GH}}{\longrightarrow}(Y,d_\infty)$ and
$(M_k,\tilde{g}_k)\stackrel{d_{GH}}{\longrightarrow}(Y,\tilde{d}_\infty)$.
Thus the convergent sequence of $(M_k,g_k)$ coincides with that of
$(M_k,\tilde{g}_k)$. That's, for $x_k\in M_k$ and $x_\infty\in Y$,
$x_k\rightarrow x_\infty$ under
$(M_k,g_k)\stackrel{d_{GH}}{\longrightarrow}(Y,d_\infty)$ iff
$x_k\rightarrow x_\infty$ under
$(M_k,\tilde{g}_k)\stackrel{d_{GH}}{\longrightarrow}(Y,\tilde{d}_\infty)$.

Next we show that $y$ is a singular point of $(Y,d_\infty)$ iff $y$
is a singular point of $(Y,d_\infty)$. This means that $\mathcal{S}$
is the singular set of $(Y,d_\infty)$. It suffices to prove the
following claim.
\begin{claim}\label{c37}
The tangent cones of $(Y,d_\infty)$ coincide with that of
$(Y,\tilde{d}_\infty)$ up to rescalings by constants.
\end{claim}
\begin{proof}[Proof of the Claim]
Let $y\in Y$ and $r_j\rightarrow0$, we want to show that, modulo a
resacling, the limits
$$\lim_{j\rightarrow\infty}(Y,r_j^{-1}d_\infty,y)=\lim_{j\rightarrow\infty}(Y,r_j^{-1}\tilde{d}_\infty,y)$$
in the pointed Gromov-Hausdorff topology, if either the limit
exists, that's, the associated tangent cone at $y$ with respect to
$d_\infty$ and $\tilde{d}_\infty$ are the same.

Passing a subsequence if necessary, the potentials $f_k$ converge to
a Lipschitz function $f_\infty$ on $Y$. Then
$\tilde{g}_\infty=e^{-\frac{2}{n-2}f_\infty}g_\infty$ on
$\mathcal{R}$. By (\ref{e35}), we have
$$Lip(f_k),Lip(f_\infty)\leq C_0$$ for some $C_0=C_0(n,D)$. Let $\rho$ be any fixed positive constant.
The distance functions $d_\infty$ and $\tilde{d}_\infty$ satisfy the
relative comparison on $B_{g_\infty}(y,r_j\rho)$,
$$\beta_j^{-1}\leq\alpha(y)\cdot\frac{\tilde{d}_\infty(x_1,x_2)}{d_\infty(x_1,x_2)}\leq
\beta_j,$$ for all $x_1,x_2\in B_{g_\infty}(y,r_j\rho)$, where
$\alpha(y)=e^{\frac{1}{n-2}f_\infty(y)}$ is a fixed constant, while
$\beta_j=e^{\frac{3}{n-2}C_0r_j\rho}$ tends to $1$ as
$j\rightarrow\infty$.

Let $(Y_y,d_y,o)=\lim_{j\rightarrow\infty}(Y,r_j^{-1}d_\infty,y)$.
By the definition of convergence, for any $\epsilon>0$, there exists
an $\epsilon$-approximation
$\psi_j:(B_{d_y}(o,\rho),d_y)\rightarrow(B_{d_\infty}(y,r_j\rho),r_j^{-1}d_\infty)$
for all $j$ large enough. Thus, whenever $j$ is large enough, the
map $\psi_j$ defines a
$10((\beta_j-1)\rho+\epsilon\beta_j)$-approximate from
$(B_{d_y}(o,\rho),d_y)$ to
$B_{\alpha(y)\tilde{d}_\infty}(y,(1+\beta_j^{1/3})r_j\rho)$. This
shows that
$$(Y,\alpha(y)r_j^{-2}\tilde{g}_\infty,y)\stackrel{d_{GH}}{\longrightarrow}(Y_y,d_y,o)$$
in the pointed Gromov-Hausdorff sense. It is equivalent to that
$$(Y,r_j^{-2}\tilde{g}_\infty,y)\stackrel{d_{GH}}{\longrightarrow}(Y_y,\alpha(y)^{-1}d_y,o)$$
in the pointed Gromov-Hausdorff sense. This proves of the claim.
\end{proof}

At last, we confirm the smoothness of $g_\infty$ on $\mathcal{R}$
and finishes the proof of Theorem \ref{t1}. We will use the
pseudolocality theorem in the argument.

\begin{claim}
$g_\infty$ is smooth and satisfies a shrinking Ricci soliton
equation on $\mathcal{R}$.
\end{claim}
\begin{proof}[Proof of the Claim]
Given any small number $r>0$, define
$$K_r=\{x\in \mathcal{R}|d_{\infty}(x,\mathcal{S})\geq r\}.$$
Then $K_r\subset K_k$ for all $k$ large enough.

Noting that $g_\infty$ is a $C^\alpha$ Riemannian metric on
$\mathcal{R}$, we have for some small constant
$\rho=\rho(r)\leq\epsilon_0r$
$$\Vol(\partial\Omega)^{n}\geq(1-\frac{1}{2}\delta_0)c_n\Vol(\Omega)^{n-1},
\hspace{0.3cm}\forall\Omega\subset
B_{g_\infty}(x_\infty,\rho),x_\infty\in K_r$$ where $c_n$ is the
Euclidean isoperimetric constant and $\epsilon_0,\delta_0$ are
constants in the Pseudolocality theorem \ref{Pseudolocality}.
Passing to the sequence $(M_k,g_k)$, since $\psi_k^*g_k$ converge to
$g_\infty$ uniformly in $C^\alpha$ on $K_{\frac{r}{2}}$, we may
assume that any domain $\Omega_k\subset
B_{g_\infty}(\psi_k(x_\infty),\frac{1}{2}\rho)$, where $x_\infty\in
K_r$, satisfies the following isoperimetric inequality
$$\Vol(\partial\Omega_k)^{n}\geq(1-\delta_0)c_n\Vol(\Omega_k)^{n-1},$$
whenever $k$ is large enough.

Let $g_k(t)$ be the Ricci flow solution with initial metric
$g_k(0)=g_k$. It's well known that $g_k(t)=(1-t)\phi_k(t)^*g_k$ for
a parameter family of diffeomorphisms $\phi_k(t)\in\Diff(M_k)$ which
is generated by the gradient field of $f_k$, cf. \cite{CK}. Because
$f_k$ is $C^1$ bounded, $\phi_k(t)(\psi_k(K_{2r}))\subset\psi_k(
K_r)$ for all $t$ in a small time interval $[0,\eta]$ where
$\eta=\eta(r)$ does not depend on $k$ whenever $k$ is large enough.
Applying the Pseudolocality theorem \ref{Pseudolocality} to points
in $\psi_k(K_{2r})$ at time $t(r)=\min(\epsilon_0r,\eta)$, we get a
uniform curvature bound
$$|Rm(g_k)|(x)\leq(1-t(r))\cdot(t(r)^{-1}+(\epsilon_0r)^{-2}),\hspace{0.5cm}\forall x\in\psi_k(K_{2r}).$$
Then Shi's gradient estimate \cite{Shi} to the Ricci flow $g_k(t)$
on the ball $B_{g_k}(x,r)$, where $x\in\psi_k(K_{3r})$, gives the
higher derivation estimate for curvature
$$|\nabla^{l}Rm(g_k)|(x)\leq C(n,l,r),\hspace{0.5cm}\forall x\in\psi_k(K_{3r}),l\geq0$$
for some constant $C$ depending only on $n,r$ and positive integer
$l$, but does not depend on specified $k$ which is large enough. It
follows from Cheeger-Gromov's compactness theorem that, passing a
subsequence once again, the metrics $\psi_k^*g_k$, modulo changes of
diffeomorphisms on subsets of $M_k$, converge in the $C^\infty$
sense to $g_\infty$ on $K_{3r}$. The arbitrariness of $r>0$ implies
that $g_\infty$ is indeed smooth on the whole $\mathcal{R}$.

In view of the soliton equation (\ref{soliton}), the gradient
estimate for curvature of $g_k$ concludes the gradient estimate for
potential functions as well:
$$|\nabla^lf_k|(x)\leq\tilde{C}(n,l,r),\hspace{0.5cm}\forall x\in\psi_k(K_{3r}),$$
where $\tilde{C}$ are constants depending only on $n,r$ and $l$. On
the other hand, one has the $C^1$ uniform bound of $f_k$ over $M_k$.
Combining these we obtain that $\psi_k^*f_k$ converge to a
$C^\alpha$ function $f_\infty$ on $Y$ which is smooth on
$\mathcal{R}$. Furthermore, by the smooth convergence, the shrinking
soliton equation
$$Ric(g_\infty)+\Hess(f_\infty)=\frac{1}{2}g_\infty$$
also holds on $\mathcal{R}$. This completes the proof of the claim.
\end{proof}

\begin{remark}
One may hope that the analogy for K\"{a}hler Ricci solitons remains
true. That's, the limit space of K\"{a}hler Ricci solitons with
positive first Chern class, under (\ref{volume noncollapsing}) and
(\ref{diameter bound}), has singular set of Hausdorff codimension at
least 4. To get this result, one needs more analysis on the tangent
cone at the singular point, cf. \cite{ChCoTi}.
\end{remark}

\subsection{Integral curvature bounds}

The following corollary has appeared in several papers, see
\cite{We}, \cite{Zhang2} and \cite{ChWa}. We reprove it here as a
corollary of our main theorem.

\begin{corollary}
Let $(M_k,g_k)$ be a sequence of closed shrinking Ricci solitons
which satisfies {\rm (\ref{volume noncollapsing}), (\ref{diameter
bound})} and
\begin{equation}\label{integral curvature}
\int_{M_k}|Rm(g_k)|^{n/2}dv_{g_k}\leq\Lambda
\end{equation}
for some $\Lambda$ independent of $k$. Then $(M_k,g_k)$ converges
along a subsequence to a compact orbifold shrinking Ricci soliton.
\end{corollary}

Here, an orbifold shrinking Ricci soliton is defined to be a smooth
Riemannian orbifold which satisfying equation (\ref{soliton}) for
some smooth function in the orbifold sense.

\begin{proof}
Let $(M,g)=(M_k,g_k)$ for a fixed $k$ and $f=f_k$ be the associated
potential function, then in local normal coordinate
$(x^1,\cdots,x^n)$, the curvature tensor of $\tilde{g}$ is given by,
cf. \cite{Be},
\begin{eqnarray}\nonumber
\tilde{R}_{ijkl}=e^{-\frac{2}{n-2}f}(R_{ijkl}+\frac{1}{(n-2)^2}g\circ((n-2)\nabla^2f+d
f\otimes d f-\frac{1}{2}|\nabla f|^2g))
\end{eqnarray}
where $\circ$ denotes the Kulkani-Nomizu product which is defined by
\begin{equation}\nonumber
(u\circ
v)_{ijkl}=u_{ik}v_{jl}+u_{jl}v_{ik}-u_{il}v_{jk}-u_{jk}v_{il}
\end{equation}
for any symmetric $(2,0)$-tensors $u$ and $v$. Substituting into the
identity $\nabla^2f=\frac{1}{2}g-Ric$ and integrating over $M$, we
get
\begin{eqnarray}
&&\int_{M}|Rm(\tilde{g})|_{\tilde{g}}^{n/2}dv_{\tilde{g}}\nonumber\\
&=&\int_M|R_{ijkl}+\frac{1}{(n-2)^2}g\circ((n-2)\nabla^2f+d
f\otimes d f-\frac{1}{2}|\nabla f|^2g)|^{n/2}dv_{g_k}\nonumber\\
&\leq&C(n)\cdot\int_M(|Rm|^{n/2}+|\nabla^2f|^{n/4}+|\nabla
f|^{n/2}+|\nabla f|^n)dv\nonumber\\
&\leq&C(n)\cdot\int_M(|Rm|^{n/2}+|\nabla f|^n+1)dv.\nonumber
\end{eqnarray}
By (\ref{e35}) and (\ref{e38}) we get a uniform bound
$C=C(n,v,D,\Lambda)$ such that
$$\int_M|Rm(\tilde{g})|_{\tilde{g}}^{n/2}dv_{\tilde{g}}\leq C.$$

By Theorem \ref{t1}, together with Theorem 2.6 of \cite{An2}, we
know that $(M_k,g_k)$ converge along a subsequence to a limit metric
space, say $(Y,g_\infty)$, with finite orbifold singularities. The
metric $g_\infty$ is globally $C^0$ and smooth on $\mathcal{R}$, the
regular part of $Y$. Furthermore, $g_\infty$ satisfies the Ricci
soliton equation (\ref{soliton}) on $\mathcal{R}$ for some function
$f_\infty$ which is globally Lipschitz and locally smooth on
$\mathcal{R}$. The remaining is to show that modulo a diffeomorphic
transformation, the lifted metric of $g_\infty$ on the resolving
domain admits a smooth extension across the singular point. The
approach is standard for Einstein manifolds by now, see \cite{An1},
\cite{BaKaNa} or \cite{Ti1} for instance; the treatment for Ricci
soliton cases are similar, we refer to \cite{CaSe} and \cite{Zh} for
details.

So we finish the proof of the theorem.
\end{proof}

\begin{remark}
In dimension four, the upper bound of $L^{n/2}$-norm of curvature
tensor in the theorem can be replaced by the bound of the second
Betti number; see \cite{Zhang2}.
\end{remark}

\section{The K\"{a}hler case}

In this section we give a proof of Theorem \ref{t2}. The proof
adopts a different method, so we put it in a new section.

\begin{proof}[Proof of Theorem \ref{t2}] By Cheeger-Gromov's compactness
theorem, it suffices to give a uniform bound for the sectional
curvature of manifolds satisfying assumptions in the theorem. We
adopt a blowing up argument to prove this. The idea comes from
\cite{RuZhZh}.

Suppose otherwise, there exists a sequence of shrinking K\"{a}hler
Ricci solitons of positive first Chern classes $(M_k,g_k)$
satisfying (\ref{volume noncollapsing}), (\ref{diameter bound}) and
(\ref{integral curvature}), and a sequence of points $p_k\in M_k$
such that
$$|Rm(g_k)|(p_k)=\sup|Rm(g_k)|\rightarrow\infty$$ as
$k\rightarrow\infty$. Then the rescaled sequence of pointed
manifolds $(M_k,Q_kg_k,p_k)$, where $Q_k=|Rm(g_k)|(p_k)$, converge
smoothly along a subsequence to a complete Ricci flat K\"{a}hler
manifold $(M_\infty,g_\infty,p_\infty)$ such that
$|Rm(g_\infty)|(p_\infty)=1$. Furthermore, $g_\infty$ has Euclidean
volume growth and integral curvature bound:
\begin{eqnarray}
&&\Vol(B_{g_\infty}(p_\infty,r))\geq\kappa
r^{2n},\hspace{0.5cm}\forall
r>0,\nonumber\\
&&\int_{M_\infty}|Rm(g_\infty)|^{n}dv_{g_\infty}\leq C,\nonumber
\end{eqnarray}
for some $\kappa>0$. Recall that $\dim_\mathbb{C}(M_\infty)=n\geq3$,
by Theorem 2 in \cite{Ti2}, one knows that $(M_\infty,g_\infty)$ is
a resolution of $\mathbb{C}^n/\Gamma$ for some finite subgroup of
$SU(n)$ which acts freely on $\mathbb{C}^n\backslash\{0\}$. $\Gamma$
is nontrivial since otherwise $M_\infty$ is simply connected at
infinity, then by Theorem 3.5 in \cite{An1}, $M_\infty$ is flat,
which contradicts with $|Rm(g_\infty)|(p_\infty)=1$. In particular,
there exists a compact subvariety in $M_\infty$, say $V$, which
represents one integral homology class in $M_\infty$. Then adopting
an argument step by step as in \cite{RuZhZh}, one gets a
contradiction. This finishes the proof of the theorem.
\end{proof}

\begin{remark}
Y.G. Zhang told the author that the treatment can also be applied to
the K\"{a}hler Einstein manifolds with non-positive Einstein
constants whose K\"{a}hler class lies into the integral
cohomological classes. For the negative K\"{a}hler Einstein case,
this is proved in \cite{Ti2}.
\end{remark}








\end{document}